\documentclass[12pt,reqno]{amsart}

\usepackage[dvips]{graphicx} 

\usepackage{
hyperref
}

\usepackage{breakurl}   

\theoremstyle{definition}

\numberwithin{equation}{section}

\renewcommand*{\MR}[1]{\href{http://www.ams.org/mathscinet-getitem?mr=#1&return=pdf}{MR #1}}
\newcommand*{\ZBL}[1]{\href{http://www.zentralblatt-math.org/zmath/en/advanced/?q=an:#1&format=complete}{Zbl #1}}

\title{Exploring Felix Klein's contested modernism}

\author[Heinig]{Peter Heinig} \address{P. Heinig, Germany}
\email{heinig@ma.tum.de}

\author[Katz]{Mikhail G. Katz}\address{M. Katz, Department of
  Mathematics, Bar Ilan University, Ramat Gan 5290002
  Israel}\email{katzmik@math.biu.ac.il}

\author[Kuhlemann]{Karl Kuhlemann}\address{K. Kuhlemann, Gottfried
  Wilhelm Leibniz University Hannover, D-30167 Hannover,
  Germany}\email{kus.kuhlemann@t-online.de}

\author[Sch\"afermeyer]{Jan Peter Sch\"afermeyer}
\address{J. P. Sch\"afermeyer, Berlin, Germany}
\email{jpschaefermeyer@gmail.com}

\author[Sherry]{David Sherry} \address{D. Sherry, Department of
  Philosophy, Northern Arizona University, Flagstaff, AZ 86011,
  US}\email{David.Sherry@nau.edu}


\begin{document}


\thispagestyle{empty}



\begin{abstract}
An alleged opposition between David Hilbert and Felix Klein as modern
\emph{vs} countermodern has been pursued by marxist historian Herbert
Mehrtens and others.  Scholars such as Epple, Grattan-Guinness, Gray,
Quinn, Rowe, and recently Siegmund-Schultze and Mazzotti have voiced a
range of opinions concerning Mehrtens' dialectical methodology.  We
explore contrasting perspectives on Klein's contested modernism, as
well as Hilbert's and Klein's views on intuition, logic, and physics.
We analyze Jeremy Gray's comment on Klein's ethnographic speculations
concerning Jewish mathematicians and find it to be untenable.  We
argue that Mehrtens was looking for countermoderns at the wrong
address.
\end{abstract}

\maketitle

\tableofcontents


\section{Introduction}

Felix Klein's stature as a modern has been subject to dispute; see
Mehrtens \cite{Me90}, Bair et al.\;\cite{18b}, Franchella \cite{Fr22},
and other commentators.  Herbert Mehrtens (1946--2021), henceforth HM,
sought to cast Klein in an unenviable role of a countermodern.  HM was
subsequently provided with a platform by Goldstein, Gray, and Ritter
to present his views in~\cite{Me96}.  HM's approach has influenced a
number of authors, including Gray \cite{Gr08}, Quinn \cite{Qu12}, and
most recently Mazzotti \cite{Ma23}, as we analyze below.

In his recent obituary for HM, Reinhard Siegmund-Schultze appears to
endorse HM's dichotomy of moderns and countermoderns (as well as to a
certain extent the assignment of Klein to the latter group).  We will
examine Siegmund-Schultze's appraisal of HM's methodology in
Section~\ref{s2}, and the views of Jeremy Gray and Frank Quinn in
Sections~\ref{s41b} and \ref{s4}, respectively.  Sections~\ref{s5}
through \ref{s33b} analyze the historical context of the 1900s in
relation to the Berlin/G\"ottingen rivalry, and its implications for
the contested issue of Klein's modernism.  Section~\ref{s8} presents
our conclusions.

\section{Siegmund-Schultze on Mehrtens and his critics}
\label{s2}

Commenting on HM's views concerning intuition and modern
mathematicians, Siegmund-Schultze mentions
\begin{enumerate}\item[]
[Mehrtens'] juxtaposition of modern mathematicians (Cantor, Hausdorff,
Hilbert, etc.) and their critics who relied and insisted on
`Anschauung' (Intuition) (Poincar\'e, Brouwer, Klein, Bieberbach).%
\footnote{Siegmund-Schultze \cite[p.\;27]{Si22}.}
\end{enumerate}
Siegmund-Schultze's comment is a reference to HM's dichotomy of
moderns versus countermoderns, with Klein listed on the countermodern
side.  Siegmund-Schultze seems to have sensed that there is a problem
with HM's application of the dichotomy to Klein:
\begin{enumerate}\item[]
[Felix] Klein stood \emph{in the middle} between proto-typical
`counter-moderns' and `moderns,' which is basically due to the fact
that Klein contributed considerably both to mathematical and technical
modernization, but, at the same time, to criticism of both as well.%
\footnote{Op. cit., p.\;28; emphasis added.}
\end{enumerate}

Did Klein stand ``in the middle'' as Siegmund-Schultze claims, or was
he a modern?  Siegmund-Schultze cites unfavorable evaluations of HM's
methodology by Rowe and Grattan-Guinness in the following terms:
\begin{enumerate}\item[]
Herbert's \emph{pioneering and novel} distinction between ``moderns''
and ``countermoderns'' has been misunderstood by some
\emph{superficial readers} of the book, and has been criticized as
misleading by some more thorough readers (Grattan-Guinness, 2009;
Rowe, 1997).%
\footnote{Ibid.; emphasis added.  Siegmund-Schultze does not specify
who such ``superficial readers'' may have been, nor in what way they
may have ``misunderstood'' Herbert's distinction.  We are familiar
with only one full-length recent critical study of HM's distinction,
namely Bair et al.\;\cite{18b}.  However, it is questionable whether
Siegmund-Schultze would have had the study \cite{18b} in mind, as some
of the arguments of \cite{18b} overlap with those of Grattan-Guinness
and Rowe, which Siegmund-Schultze apparently finds to be ``more
thorough'' than those of the unnamed superficial readers.}
\end{enumerate}
Siegmund-Schultze's reference is to the criticisms by Grattan-Guinness
in~\cite{Gr09} and Rowe in~\cite{Ro97}.  Indeed, neither of the pair
would have necessarily endorsed a description of HM's distinction as
``pioneering and novel,'' particularly in application to Felix Klein,
as we analyze in Section~\ref{s21}.

\subsection{Shocked reaction}
\label{s21b}

Siegmund-Schultze's attempt to put a good face on HM's distinction
becomes all the more surprising when compared with his initial
reaction to HM's book, which was that of shock:
\begin{enumerate}\item[]
It came as a shock to me since it was certainly not aiming in the
direction of our project.%
\footnote{Siegmund-Schultze \cite[p.\;27]{Si22}.}
\end{enumerate}
He and HM were working on a joint project in the 1980s that was
supposed to result in a book, when Siegmund-Schultze received a draft
copy of HM's book in june 1989 and was apparently shocked by its
marxist and semiotic content.  Siegmund-Schultze then quotes HM and
concludes: ``the quotation describes quite clearly Herbert’s
individual scientific development and gives the essential explanation
why \emph{our project had to fail}.''%
\footnote{Siegmund-Schultze \cite[p.\;30]{Si22}; emphasis added.}

On the subject of HM's dichotomy, Grattan-Guinness wrote, rather
critically:
\begin{enumerate}\item[]
On the place of modernism in mathematics as a whole, the claim [by
  Mehrtens] that `Modernism and Antimodernism make between them the
success history of mathematics in the twentieth century' [M,
  dustjacket]%
\footnote{This bracketed comment is in Grattan-Guinness' text.}
requires us to rejudge as failures or footnotes the developments of,
for example, topology, mathematical statistics, and numerical methods!%
\footnote{Grattan-Guinness \cite[pp.\;9--10]{Gr09}.}
\end{enumerate}

\subsection{Rowe on Mehrtens' dialectics}
\label{s21}

Rowe points out that HM's book ``was written during the midst of the
\emph{Historikerstreit} of the late 1980s when debates regarding the
historical significance of the Holocaust were in full flame.''%
\footnote{Rowe \cite[p.\;517]{Ro13}.}
He notes that some of HM's protagonists can be viewed as modern in one
sense and counter-modern in another, and that
\begin{enumerate}\item[]
Klein and Hilbert shared strikingly similar views about what
constituted good and bad mathematics.%
\footnote{Rowe \cite[p.\;536]{Ro97}.}
\end{enumerate}
Rowe also questions HM's claim concerning antisemitism.  HM claimed
that the movement that embraced \emph{Anschauung} in the tradition of
Klein germinated into the virulent form of antimodernism and
anti-Semitism.%
\footnote{For a recent reaction by Mazzotti to such a claim, see the
main text at note~\ref{f19}.}
Writes Rowe:
\begin{enumerate}\item[]
Far more provocative and central for Mehrtens's whole approach to
countermodernism are his claims regarding the impact Klein's views
exerted on subsequent events. According to Mehrtens, the movement that
embraced Anschauung in the tradition of Klein--a movement that gained
momentum after 1905 with the adoption of the Meran proposals that
called for the incorporation of anschauliche methods in mathematics
instruction--germinated into the virulent form of antimodernism and
anti-Semitism espoused by Ludwig Bieberbach, the leading proponent of
Aryan mathematics.%
\footnote{Op. cit., p.\;5.}
\end{enumerate}
Rowe's verdict is unequivocal:
\begin{enumerate}\item[]
Unfortunately, Mehrtens tells us only how this transition \emph{might
have} happened instead of documenting that it \emph{actually did}
happen.%
\footnote{\label{f8}Ibid.; emphasis added.}
\end{enumerate}
Given the weakness of HM's positing such a link, one can well wonder
what the source of HM's animus toward Klein might have been.  Rowe
provides a hint when he speaks of ``\emph{dialectical} elements so
prominent in Mehrtens’s interpretation''%
\footnote{Rowe \cite[p.\;519]{Ro13}; emphasis added.}
and mentions that HM
\begin{enumerate}\item[]
paints a harsh picture that underscores moral bankruptcy and
dehumanizing influences, as when mathematicians willingly lent their
technical expertise to the military.%
\footnote{Rowe \cite[pp.\;516--517]{Ro13}.}
\end{enumerate}
The academic-industrial complex was the subject of a 1908 cartoon
reproduced in HM's book.%
\footnote{For a discussion of the cartoon see Bair et
al.~\cite[Section~4.7]{18b}.}
%
%
HM's distaste for the Academy's cooperation with the industry and the
military arguably stemmed from HM's political commitments and colored
his evaluation of Klein's relation to modernism.  As noted by Gray,
\begin{enumerate}\item[]
Mehrtens's critique was written in a post-Marxist spirit, influenced
by such writers as Foucault.%
\footnote{Gray \cite[p.\,10]{Gr08}.}
\end{enumerate}

Rowe adds that ``Mehrtens evinces little interest in such `details'
that would clutter up the main picture he wishes to convey.''%
\footnote{Rowe \cite [p.\;4]{Ro97}.}
When HM's crude marxism and semiotic obfuscation escalate, Rowe does
not shy away from calling him to task:
\begin{enumerate}\item[]
Mehrtens goes off the deep end, for example, when he tries to find
parallels between the formal languages of modern mathematics and the
commando rhetoric of contemporary dictatorships.%
\footnote{Rowe \cite[p.\;520]{Ro13}.}
\end{enumerate}

\subsection{Mehrtens as seen by Epple, Tobies and Mazzotti}

Epple is similarly critical of HM's insufficient clarity in the
definition of the terms \emph{modern} and \emph{countermodern}:
\begin{enumerate}\item[]
Mehrtens' book is an example of a very elaborated kind of external
historiography.  His sources are mainly the programmatic declarations
of the mathematicians involved and the documents of their
institutional activities.  Mehrtens does not attempt to analyze some
of the more advanced productions of modernist or counter-modernist
mathematicians {\ldots}\;Thus we are left in a somewhat unclear
position if we accept his narrative.%
\footnote{Epple \cite[p.\,191]{Ep97}.}
\end{enumerate}
Epple goes on to propose his own interpretation of what the dichotomy
of modern vs countermodern could mean.  We will not deal with Epple's
definition because our main concern is HM's distinction rather than
Epple's.  HM's definition is characteristically slippery, but at least
two features of mathematical modernism seem indisputable:

\begin{enumerate}
\item
a tendency to unify, rather than segregate, fields within mathematics;
and
\item
in seeking to understand what constitutes mathematical modernism, we
cannot ignore what leading 20th century mathematicians thought
characterized modern mathematics.%
\footnote{\label{f13}More on this in the main text at note~\ref{f47}.}
\end{enumerate}
By the unifying criterion, Klein wins hands down, because of his
efforts to unify pure and applied mathematics, as detailed by Tobies.%
\footnote{Tobies \cite[Chapters 7.7 and 8.1]{To21}.}

In a similar vein, Tobies writes:
\begin{enumerate}\item[]
Contrary to Herbert Mehrtens's (1990) classification of mathematicians
as either modern (Hilbert et al.) or anti-modern (Klein et al.), it
seems to me as though, in Klein's case, it is possible to speak of a
particular sort of modernity on whose basis new domains such as modern
numerical analysis, actuarial mathematics, and financial mathematics
were able to develop.%
\footnote{Tobies \cite[p.\,10]{To21}.}
\end{enumerate}
Thus, Tobies tends to assign Klein to the \emph{moderns}.  Meanwhile,
Mazzotti writes:
\begin{enumerate}\item[]
[I]n one of the \emph{strongest} claims of the book, [Mehrtens] argues
for an essential continuity between this [i.e., Felix Klein's] strain
of Wilhelmian mathematics and the anti-modern and anti-Semitic
campaign of the proponents of Aryan mathematics under the Nazi
regime.%
\footnote{\label{f19}Mazzotti \cite[p.\;277, n.\;20]{Ma23}; emphasis
added.}
\end{enumerate}
As noted by Rowe, HM did not \emph{argue for} but rather merely
\emph{postulated} such an ``essential continuity.''  With regard to
the qualifier ``strongest'' it is unclear whether it is being used in
the sense of \emph{outrageous} or in the sense of \emph{compelling} by
Mazzotti, who does not conceal MH's influence on his work:
\begin{enumerate}\item[]
In this book I have clearly followed Mehrtens's invitation to explore
the politics of mathematical modernity, {\ldots}%
\footnote{Mazzotti \cite[p.\;238]{Ma23}.  Mazzotti describes Mehrtens'
\cite{Me90} as ``an ambitious and original book {\ldots}\;that deploys
high-powered semiotic and philosophical tools to study the discourse
of the autonomy and meaning of mathematics'' \cite[p.\;237]{Ma23}.}
\end{enumerate}

In contrast to Grattan-Guinness, Rowe, and Tobies, Jeremy Gray adopts
a more sympathetic view toward HM's approach to Felix Klein.  In a
2008 book, Gray voices criticisms of Klein that echo HM's.  We will
examine Gray's portrayal of Klein in Section~\ref{s41b}.

\section{Gray on Klein's `notorious' paragraph}
\label{s41b}

Like HM, Jeremy Gray attempts to tackle modernism in his 2008 book
\cite{Gr08}.  Grattan-Guinness expresses scepticism with regard to
Gray's definition of modernism in the following terms:
\begin{enumerate}\item[]
[W]hile Gray relates his definition to a text by Guillaume Apollinaire
on cubist art, it is surely incomplete and so too wide; for example,
there is no reference to mathematics at all, and no time dependence.%
\footnote{Grattan-Guinness \cite[p.\;4]{Gr09}.}
\end{enumerate}

HM wrote that
\begin{enumerate}\item[]
Theodor Vahlen
{\dots}~cited \emph{Klein's racist distinctions} within mathematics,
and sharpened them into open antisemitism.%
\footnote{Mehrtens \cite[pp.\;310--311]{Me90}; translation ours;
emphasis added.}
\end{enumerate}
Even though HM is only mentioned once in a footnote in Gray's
Section~4.2.1.1 on pages 197--199, we will see that HM's methodology
exerted significant influence on Gray's portrayal of Klein.  Most
commentators today would agree that Klein's distinctions, as discussed
below, were \emph{racial} rather than \emph{racist}.  The implied
connection between Klein's offhand comments on distinct ethnic styles
in mathematics, on the one hand, and anti-Semitism, on the other,
bears closer scrutiny.  In this section, we will analyze Jeremy Gray's
take on Klein's comments.

\subsection{`Teutons' and Jews}

As is well known, Klein was instrumental in securing positions for a
number of qualified Jewish mathematicians in German universities.%
\footnote{See e.g., Rowe \cite{Ro86}.}
Regrettably, Gray tends to downplay Klein's valiant efforts on behalf
of Jewish mathematicians.  One would suspect little of such efforts by
reading Gray:
\begin{enumerate}\item[]
Then came a paragraph that became notorious.$^{48}$\, {\small
    [here 48 is a footnote in Gray's book; see Section~\ref{s32}]}
  Klein speculated that the exactness of spatial intuition might vary
  between individuals along racial lines.  \emph{Teutons} seemed to
  have strong, naive spatial intuition, while the critical, purely
  logical sense was more fully developed in the Latin and Hebrew
  races.%
\footnote{Gray \cite[p.\;198]{Gr08}; emphasis added.}
\end{enumerate}
It is worth analyzing the accuracy of Gray's passage in detail.  In
fact, there are some significant inaccuracies in his narrative.  At
variance with Gray's attempted summary, Klein himself never spoke of
\emph{Teutons}, but only of \emph{the Teutonic race}.  Klein's passage
reads as follows:
\begin{enumerate}\item[]
Finally, it must be said that the degree of exactness of the intuition
of space may be different in different individuals, perhaps even in
different races.  It would seem as if a strong naive space-intuition
were an attribute pre-eminently of the Teutonic race, while the
critical, purely logical sense is more fully developed in the Latin
and Hebrew races.  A full investigation of this subject, somewhat on
the lines suggested by Francis Galton in his researches on heredity,
might be interesting.%
\footnote{Klein \cite[pp.\;45--46]{Kl94}.}
\end{enumerate}
In Klein's case, there was no difference between words and actions:
whatever ethnographic comments he may have made (and whatever their
depth or lack thereof), the comments were all positive rather than
critical; namely, he pointed out the positive aspects of having
distinct styles in mathematics, as in the case of his comments on
Sylvester.%
\footnote{See Bair et al.~\cite[Section 2.2]{18b}.}
His words were consistent with his actions on behalf of many Jewish
mathematicians, as noted by Tobies (see Section~\ref{s32}).

In the passage quoted above, Klein is not using the term
\emph{Teutons}, which would be the name of a \emph{people} but rather,
Klein refers to \emph{the Teutonic race} which indicates that Klein is
simply dabbling in the ethnographic theorizing that was common for
this historical period.  Gray continues:

\begin{enumerate}\item[]
A fuller investigation along the lines of Galton's study of heredity
might be interesting, he [i.e., Klein] concluded.  Given Klein's
evident preference for intuition, it is \emph{hard not to read} this
passage as being somewhat anti-Semitic.%
\footnote{Gray \cite[p.\;198]{Gr08}; emphasis added.}
\end{enumerate}
Gray finds that it is ``hard not to read'' Klein's paragraph as
``being somewhat anti-Semitic.''  Apparently some find it \emph{easy
to read} the passage as anti-Semitic.  We find it hard not to ask,
\emph{who} exactly would find it easy to read Klein's paragraph in
this light?  Certainly not the historian Segal \cite{Se03}.%
\footnote{Regarding Segal's position, see further in Bair et
al.~\cite{18b}.}
While Ludwig Bieberbach (who appears to have truly believed in his
crude aryanism) would have surely found it easy to read Klein's
passage in that light, we reject the attribution of indirect
culpability to Klein by Gray following HM.  The idea that a modern
historian would find such an association of ideas hard to resist
strikes us as peculiar.  Gray continues:
\begin{enumerate}\item[]
[1] Moreover, Kronecker, Klein's examplar of the purely logical
mathematician, was Jewish; Peano, a Latin, had already been mentioned
by Klein as an axiomatic mathematician.  [2]~Forty years later, these
remarks were brandished by the Nazi Bieberbach as `proof' that `good'
German mathematicians such as Klein and Bieberbach were racially
superior.%
\footnote{Gray \cite[p.\;198]{Gr08}; numerals [1] and [2] added.}
\end{enumerate}
As per Gray's sentence [1], Klein paints Kronecker and Peano as
logicians and axiomatizers; as per sentence [2], Klein's remarks
allegedly provided grist for the mill of Bieberbach's racist
propaganda.  Does Gray take pains to distance himself from the
inference that Klein could bear partial responsibility for enabling
such propaganda?  On the contrary: Gray seems to come dangerously
close to falling into the same trap as HM in confusing possibility and
actuality (as noted by Rowe):%
\footnote{See the main text at note~\ref{f8}.}
\begin{enumerate}\item[]
Though Bieberbach made poisonous nonsense of it, there is no denying
that Klein's remarks partake of the \emph{nineteenth-century glibness}
with which whole races are pigeon-holed.  Psychology has ever lent
itself to the representation of contemporary prejudice as the latest
fruit of science.%
\footnote{Ibid.; emphasis added.}
\end{enumerate}
There may be no denying the glibness, but one can certainly deny that
Klein's passage was ``somewhat anti-Semitic,'' contrary to Gray's glib
claim.  The glibness with which a 21-century commentator dismisses an
entire academic discipline (namely, \emph{Psychology}) as lending
itself ``to the representation of contemporary prejudices'' is
surprising.

\subsection{A belated footnote on academic politics}
\label{s32}

In footnote~48 on the same page, Gray tries to hedge his bets with
regard to Klein:
\begin{enumerate}\item[]
$^{48}$It is discussed in Rowe 1986%
\footnote{This is a reference to Rowe \cite{Ro86}.}
and Mehrtens 1990, 215--218, who also note that elsewhere in his
\emph{Entwicklung} Klein is clear that the emancipation of the Jews in
Germany in 1812 was quickly fruitful for mathematics, and that Klein
was not anti-Semitic in his \emph{academic politics}, rather the
opposite.%
\footnote{Gray \cite[p.\;198]{Gr08}; emphasis added.}
\end{enumerate}
Gray's apparent conclusion is that it is only in his ``academic
politics'' that Klein wasn't antisemitic, retaining the unstated
possibility of applying the slur in other areas.  Regrettably, Gray
fails to distance himself sufficiently from HM's comments on alleged
affinities between Bieberbach's views on \emph{Jewish mathematics} --
and those of Klein.  Writes Tobies:
\begin{enumerate}\item[]
Klein’s classification here was not intended to denigrate any
mathematical approach or ``race.''  He assigned equal value to the
different ways of achieving knowledge.  Nevertheless, after his death,
statements of this sort were politically appropriated and misconstrued
for anti-French or anti-Semitic purposes.%
\footnote{Tobies \cite[p.\;491]{To21}.}
\end{enumerate}
Tobies' appraisal appears the more balanced.

Gray claims further that ``Klein's own mathematics always inclined to
the sort of imprecision [Oskar] Perron was adamantly against''%
\footnote{Gray \cite[p.\;277]{Gr08}.}
but gives no evidence to back up his claim, and instead hedges his
bets in the following terms:
\begin{enumerate}\item[]
Even if Perron had not had Klein in mind [in Perron's 1911 lecture],
Klein's example was the kind of thing he was against%
\footnote{Ibid.}
\end{enumerate}
but again without presenting any evidence.  On the contrary, we find
Perron's view concerning the limited usefulness of geometric
intuition%
\footnote{See Gray \cite[p.\;275]{Gr08}.}
in agreement with Klein's view on this topic, as laid out in his
Evanston lecture on space-intuition in 1893, where he warned of the
pitfalls of naive intuition which can be unreliable.%
\footnote{See Klein \cite[p.\;225]{Kl93}.  For a detailed analysis,
see Franchella \cite[pp.\,1061--1062]{Fr22}.}
See also Section~\ref{s33b}.

\section{Quinn's portrayal of Klein}
\label{s4}

Another author who appears to have been influenced by HM's
interpretation, via Gray, in his thinking about Klein is Frank Quinn.
In this section we examine the image of Klein as allegedly a
traditionalist who opposed the new methods in mathematics and clung to
19th century values, that emerges from Quinn's writings.  Quinn's more
specific target is Klein's book \emph{Elementarmathematik vom
h\"oheren Standpunkte aus} \cite{Kl08}.  Quinn's article \cite{Qu12}
did not go unnoticed: it was selected for the collection \emph{The
Best Writing on Mathematics 2013} published by Princeton University
Press.%
\footnote{See Pitici \cite[pp.\;175--190]{Pi14}.}
We will argue that it was, alas, \emph{not the best}.  Given below is
a sample of Frank Quinn's comments on Klein.

\begin{enumerate}
\item
``[D]uring the transition [period] some traditional mathematicians,
  most notably Felix Klein, were very influential in education and
  strongly imprinted nineteenth-century values on early
  twentieth-century education reforms\ldots{}''%
\footnote{Quinn \cite[section\;1.5, p.\;4]{Qu11}.}
\item
``[Felix Klein's] 1908 book \emph{Elementary Mathematics from an
Advanced Viewpoint} was a virtuoso example of
  \emph{nine\-teenth-century methods} and did a lot to cement their
  place in education. The `Klein project' {\ldots}~is a contemporary
  international effort to update the topics in Klein's book but has no
  plan to update the methodology.  In brief, traditionalists lost the
  battle in the professional community but won in education.''%
\footnote{Quinn \cite[p.\;33]{Qu12}; emphasis added.}
\item
``Felix Klein was still denouncing the new methods in the 1920s, and
because his views were not only unrefuted but almost unchallenged,
outsiders accepted them as fact.''%
\footnote{Ibid., p.\;36.}
\item
``Traditionalists like Poincar\'e and Felix Klein {\ldots}~objected
  because they wanted to use physically or experientially based
  methodology.''%
\footnote{Quinn  \cite[section\;5.2.3, p.\;39]{Qu11}.}
\item
``Klein was not at all bothered by non-Euclidean geometry. However
  this may have been because he saw a place for them in a deeper
  understanding of the physical world; he remained committed to the
  physical-realist philosophy.  The new [i.e., Hilbert's] view is that
  mathematics is a world of its own, with wonderful things not dreamt
  of in their philosophy.''%
\footnote{Quinn \cite[Section\;8.1.1, p.\;47]{Qu11}.  Here Quinn both
  misrepresents Klein's views and distorts Hilbert's by implying that
  the latter would subscribe to a view of mathematics as ``a world of
  its own;'' see Section~\ref{s33b}.}
\end{enumerate}
We will summarize Quinn's position in Section~\ref{s42}.

\subsection{Summary of Quinn's position}
\label{s42}

As sourced above, Frank Quinn claims that Felix Klein (i) practiced
19th century methods as opposed to Quinn's rigorous Core Mathematics,%
\footnote{Quinn's Core Mathematics involves a combination of
Weierstrassian \emph{Epsilontik} and a focus on precise
\emph{definitions} in an explicit axiomatic framework.  Quinn
champions the late 1800s as the formative period of Core Mathematics,
mentioning the work of Weierstrass, Peano, Frege, and Hilbert.}
%
(ii) opposed Weierstrassian methods, and (iii) wanted to use a
physically based methodology.  As detailed in \cite{18b}, the truth
about Klein is precisely the opposite: he \emph{rejected} 19th century
methods, welcomed the Weierstrassian revolution,%
\footnote{Weierstrass (1815--1897) of course lived in the 19th
  century, but the Weierstrassian revolution in rigor in the
  foundations of analysis is universally considered part of the
  modernist transformation of mathematics in the 20th century.}
and clearly understood that mathematics deals only with
\emph{representations} of reality.  Writes Corry:
\begin{enumerate}\item[]
[I]n 1872, Felix Klein--using ideas developed by Cayley and aware of
the independence of the parallel axiom in Euclidean
geometry--reformulated the whole issue in terms of axiomatic
interdependence: projective geometry should be axiomatically defined,
so that the various metrical geometries can be derived by addition of
new axioms.%
\footnote{Corry \cite[p.\;155]{Co04}.}
\end{enumerate}

\subsection{Quinn's \emph{footwork} claim}
\label{s43}

Some of Quinn's most puzzling claims concern Klein's students:
\begin{enumerate}
\item[(a)]
Early in the twentieth century Klein was no longer a leader in
mathematics because no one would follow him: talented young people saw
giving up physical intuition in exchange for technical power as a good
deal and \emph{voted with their feet}.  It was then that he turned his
focus to education\ldots%
\footnote{Quinn \cite[Section\;8.1.3, p.\;47]{Qu11}; emphasis added.}
\item[(b)] Traditionalists, including Felix Klein {\ldots}~had no
  audience in the professional community because young people
  \emph{voted with their feet} for the new ways.%
\footnote{Quinn \cite[Section\;1.5, p.\;9]{Qu11b}; emphasis added.}
\end{enumerate}
How would one evaluate Quinn's \emph{footwork} 2011 claim in items (a)
and~(b)?  The useful \emph{Mathematics Genealogy Project} is in the
public record, and records no fewer than nine students who defended
their PhD under Klein in the first decade of the 20th century.  They
are Carl Wieghardt, Conrad M\"uller, Max Winkelmann, Anton Timpe,
W.\;Ihlenburg, Julio Rey Pastor, Hermann Rothe, Ludwig Bieberbach, and
Erwin Freundlich.%
\footnote{See
  \url{https://www.genealogy.math.ndsu.nodak.edu/id.php?id=7401}}
Klein's more than 77500 mathematical descendants (as of january 2024)%
\footnote{Of this number, over 41000 are through Hilbert, who was
Klein's mathematical grandson through Lindemann.}
apparently have not noticed the general stampede away from Klein
implied by Quinn.

Moreover, Hermann Weyl, who obtained his doctorate under David
Hilbert, received detailed assistance and mathematical advice from
Klein in writing his distinctly modern book on Riemann surfaces
\cite{We13}.%
\footnote{See especially pages V and IX there and
  \cite[Section\;1.10]{18b} for additional details.}

Quinn's \emph{footwork} claim dovetails with Quinn's \emph{tale}, but
such footwork does not seem to be part of the historical record.
Quinn's claim concerning young people allegedly walking away from
Felix Klein appears to be contrary to historical evidence.

Klein indeed emphasized the importance of intuition, particularly in
education, as he will in more detail in his book \cite{Kl08}.
However, Klein's stance on education is similar to that of many a
modern mathematician who rejects Quinn's ideology on mathematics
education.

Thus, Saunders Mac Lane wrote: 
\begin{enumerate}\item[]
The sequence for the understanding of mathematics may be:
\emph{intuition, trial, error, speculation, conjecture, proof}.%
\footnote{Mac Lane in \cite[p.\;191]{At94}; emphasis in the original.}
\end{enumerate}
Dismissing Klein's take on education as a throwback to the 19th
century, as Quinn does, is as unfounded as dismissing Michael Atiyah's%
\footnote{Atiyah et al.~\cite{At94}.}
on the same grounds.  In a related vein, Terry Tao distinguishes three
stages in a student's mathematical education: (1)\;pre-rigorous stage;
(2) rigorous stage; (3) post-rigorous stage.  Tao describes the
pre-rigorous stage as follows:
\begin{enumerate}\item[]
In [the pre-rigorous stage] mathematics is taught in an informal,
intuitive manner, based on examples, fuzzy notions, and hand-waving.
(For instance, calculus is usually first introduced in terms of
slopes, areas, rates of change, and so forth.) The emphasis is more on
computation than on theory. This stage generally lasts until the early
undergraduate years.%
\footnote{Tao \cite{Ta09}.}
\end{enumerate}
By Quinn's lights, Tao's first stage would apparently amount to a
counterproductive throwback to the 19th century.%
\footnote{\label{f47}The pertinence of taking into account the views
of leading 20th century mathematicians was already mentioned in the
main text at note~\ref{f13}.}

\section{Hilbert and Klein: allies or adversaries?}
\label{s5}

While it is true that Klein defended an older mathematical tradition
than the one envisioned by Hilbert, the view of Hilbert and Klein as
adversaries in an alleged struggle of \emph{modern} versus
\emph{countermodern}, found in the writings of Mehrtens \cite{Me90},
Gray \cite{Gr08}, and Quinn \cite{Qu12} as analyzed in
Sections~\ref{s2}, \ref{s41b}, and \ref{s4}, can be challenged.
Klein's determined campaign to bring Hilbert to G\"ottingen finally
paid off in 1895:
\begin{enumerate}\item[]
Klein clearly knew what he wanted as well as what he was getting, and,
for his part, Hilbert gladly joined forces with him in his
long-standing battle with the Berliners, a fight that continued to
rage during the post-Kroneckerean era.%
\footnote{Rowe \cite[pp.\;10--11]{Ro13}.}
\end{enumerate}

\subsection{United against the Berlin school}

Rowe notes that Klein and Hilbert were \emph{allies} in resisting a
``purist trend toward specialization'' emanating from the Berlin
school, dominated by Weierstrass, Kummer, and Kronecker:
\begin{enumerate}\item[]
[Klein and Hilbert] forged an intellectual partnership that not only
overturned the balance of power within German mathematics but also
significantly altered its research orientation and the relationship
between mathematics and other scientific and technological
disciplines.%
\footnote{Rowe \cite[p.\;186]{Ro89}.}
\end{enumerate}
The Berlin school was characterized by a ``belief that arithmetic was
the true font of all mathematical knowledge and that other branches of
mathematics were to be regarded as subordinate to it.''%
\footnote{Rowe \cite[p.\;187]{Ro89}.}
Indeed,
\begin{enumerate}\item[]
It was this \emph{purist trend toward specialization} in mathematical
research that Klein and Hilbert resisted, turning instead to a
broader, all-embracing approach that stressed the vitality of unifying
ideas {\ldots}\;Both were universalists with grandiose visions for the
mathematics of the future, visions they felt were threatened by
institutional and intellectual forces alike. And although their
aspirations for mathematics diverged in certain important respects,
underlying them was a fundamental agreement on most key issues.%
\footnote{Ibid.; emphasis added.}
\end{enumerate}
Echoes of such Berlin--G\"ottingen tensions can be detected in
Hilbert's 1900 Paris lecture (its idea is traceable to an invitation
by Hermann Minkowski):
\begin{enumerate}\item[]
While insisting on rigour in the proof as a requirement for a perfect
solution of a problem, I should like, on the other hand, to oppose the
opinion that only the concepts of analysis, or even those of
arithmetic alone, are susceptible of a fully rigorous treatment.  This
opinion, occasionally advocated by \emph{eminent men}, I consider
entirely erroneous.%
\footnote{Hilbert \cite[p.\;442]{Hi00}; emphasis added.}
\end{enumerate}
The ``eminent men'' here is a thinly veiled allusion to the Berliners.
Having spelled out the error of the analysis \emph{\"uber alles} view
as he saw it, Hilbert continued to emphasize the vital link to
geometry and mathematical physics:
\begin{enumerate}\item[]
Such a one-sided interpretation of the requirement of rigour would
soon lead to the ignoring of all concepts arising from geometry,
mechanics and physics, to a stoppage of the flow of new material from
the outside world, and finally, indeed, as a last consequence, to the
rejection of the ideas of the continuum and of the irrational number.
But what an important nerve, vital to mathematical science, would be
cut by the extirpation of geometry and mathematical physics!  On the
contrary I think that wherever, from the side of the theory of
knowledge or in geometry, or from the theories of natural or physical
science, mathematical ideas come up, the problem arises for
mathematical science to investigate the principles underlying these
ideas and so to establish them upon a simple and complete system of
axioms, that the exactness of the new ideas and their applicability to
deduction shall be in no respect inferior to those of the old
arithmetical concepts.%
\footnote{Ibid.}
\end{enumerate}
The ``eminent men" criticized by Hilbert were Weierstrass, Kummer,
Kronecker (all three were dead by the time of Hilbert's Paris
lecture), Frobenius, Schwarz, and others of the Berlin school.  Writes
Rowe:
\begin{enumerate}\item[]
The Berlin mathematicians had long treated geometrical arguments as
foreign to the rigorous standards demanded by algebra and analysis
{\ldots}\;Diagrams and other aids to the imagination were considered
taboo.  Hilbert's aim was to dispel this orthodoxy by showing that the
foundations of geometry were every bit as rigorous as those based on
the properties of number systems.%
\footnote{Rowe \cite[p.\,14]{Ro13}.}
\end{enumerate}
As we will see in Section~\ref{s62}, Klein was passed over in favor of
Schwarz for promotion at Berlin.  We will compare their attitudes
toward set theory in Section~\ref{s52}.

\subsection{Fraenkel on Schwarz on set theory}
\label{s52}

In the area of set theory, we have a memorable recollection of Abraham
Fraenkel (of ZFC fame and impeccable \emph{modernist} credentials)
concerning a 1912 seminar at Berlin where Fraenkel (as a student) gave
a lecture on set theory.  Fraenkel documents Hermann Schwarz's
reaction in the following terms:
\begin{enumerate}\item[]
Schwarz took the most convenient route.  He gave interested
students complete freedom to report on whatever subjects they chose
during the seminar.  In unconscious anticipation of my future
developments, I chose as my subject the elements of set theory, which
was virtually unknown in Berlin.  My two-hour lecture stimulated
lively participation among the students.  The professor then announced
in his closing remarks that he too had heard about Georg Cantor's
questionable theories, but felt he had to warn young students about
them.%
\footnote{Fraenkel \cite[pp.\;94--95]{Fr16}.}
\end{enumerate}
Note that such memorable references to ``Georg Cantor's questionable
theories,'' assorted with the need to ``warn young students about
them,'' were uttered by Schwarz at Berlin, not Klein at G\"ottingen.
As we know, Klein fully endorsed Cantor's theory; see
Section~\ref{s63}.

\subsection{Fraenkel, Lie, Noether, Zermelo}

Fraenkel offered the following remarkable recollection on Klein:
\begin{enumerate}\item[]
My professional career began in march 1919 with {\ldots}~an invitation
to G\"ottingen to Privy Counselor Felix Klein, almost 70 years of age,
but still active as the ``foreign minister'' of German mathematics.%
\footnote{Fraenkel \cite[p.\,127]{Fr16}.}
\end{enumerate}
We may therefore add to the list of Klein's credits, the fact of
having helped launch Fraenkel's professional career.  This is a
remarkable accomplishment to Klein's credit.  Klein invited Fraenkel
to contribute an article on \emph{Number} for a collection on Gauss
that Klein was publishing.  Since Fraenkel was a rising star in the
new field of axiomatic set theory, the fact that Klein asked Fraenkel
to write an article on \emph{Number} suggests that Klein still adhered
to his view concerning the importance of the \emph{arithmetisation of
analysis}, a term he coined in his 1895 address; see \cite{Kl96}.
Possibly Klein invited Fraenkel to write an article on \emph{Number}
because Klein considered set theory to be the foundation of such an
arithmetisation--a markedly modernist stance.

Klein was instrumental in securing Sophus Lie's appointment at Leipzig
(against Weierstrass' opposition; see \cite{18b}).  We have a detailed
testimony on how the moderns Klein and Lie reacted to Weierstrass'
lectures in 1869.  Wrote Klein:
\begin{enumerate}\item[]
As far as I remember --- I arrived at Berlin in 1869 and stayed there
from 1869 to 1870 --- Weierstrass had the status of an absolute
authority in whose teachings the listeners acquiesced, [taking them]
as an indisputable norm, oftentimes without having properly
comprehended them in a deeper sense.  Doubts weren't allowed to arise,
and checking was difficult, already because of Weierstrass' giving
extraordinarily few literature references.  In his lectures he had set
himself the goal of presenting a system of neatly ordered and coherent
thoughts.  Thus he started with methodical construction, from the
ground up, and, in pursuit of his ideal of not leaving out anything,
he set out to proceed in such a way as to later only ever have to
refer to himself.%
\footnote{In the original: ``Nach meinen Erinnerungen --- ich kam 1869
nach Berlin und war 1869/70 dort --- war Weierstra{\ss}' Stellung die
einer absoluten Autorit{\"a}t, deren Lehren die Zuh{\"o}rer hinnahmen
als unanfechtbare Norm, oft ohne sie im tieferen Sinn recht
aufgefa{\ss}t zu haben. Ein Zweifel durfte nicht aufkommen, eine
Kontrolle war schon deshalb schwer m{\"o}glich, da Weierstra{\ss}
au{\ss}erordentlich wenig zitierte. Er hatte es sich in seinen
Vorlesungen zum Ziel gesetzt, ein System wohlgeordneter Gedanken im
Zusammenhang vorzutragen. So begann er mit einem methodischen Aufbau
von unten herauf und, seinem Ideal der L{\"u}ckenlosigkeit
nachstrebend, richtete er den Gang so ein, da{\ss} er in der Folge nur
auf sich selbst zur{\"u}ckgreifen brauchte.''  Klein
\cite[p.\;284]{Kl26}.}
\end{enumerate}
Klein continued:
\begin{enumerate}\item[]
Back then, I myself---now I regret it---did the same as did [Sophus]
Lie: out of contrariness, I attended none of Weierstrass' courses, but
rather always fought for my own ideas in [Weierstrass'] seminar. But
back then I did hand-copy a transcript of a course of Weierstrass' on
elliptic functions, and I used it years later when working on this
subject.

Gradually Weierstrass gained a reputation all over the scientific
world as an incomparable authority (cf.\;Mittag-Leffler,
[International] Congress [of Mathematicians in] Paris [1900],
page\;131, where Hermite writes: `Weierstrass is the master of us
all').

And yet, Weierstrass, too, was finally not spared the disappointment
of having to watch his doctrines being challenged (cf.\;letter to
[Sofia] Kovalevskaya dated 24 March 1885, communicated by
Mittag-Leffler, Acta Math., Volume 39, p.\;194ff.).%
\footnote{``Ich selbst habe damals---jetzt bedaure ich es---ebenso wie
Lie, aus Widerspruchsgeist kein Kolleg bei Weierstra{\ss} geh{\"o}rt,
sondern im Seminar immer nur eigene Gedanken verfochten. Aber eine
Vorlesung von Weierstra{\ss} {\"u}ber elliptische Funktionen habe ich
mir damals abgeschrieben und Jahre sp{\"a}ter bei meinen Arbeiten
{\"u}ber diesen Gegenstand oft benutzt.

Allm{\"a}hlich gewann Weierstra{\ss} Geltung in der gesamten
wissenschaftlichen Welt als unvergleichliche Autorit{\"a}t (vgl. dazu
Mittag-Leffler, Pariser Kongre{\ss} 1900, S.\;131, wo Hermite sagt:
`Weierstra{\ss} est notre ma{\^i}tre {\`a} tous').

Und dennoch ist zuletzt auch Weierstra{\ss} die Entt{\"a}uschung nicht
erspart geblieben, da{\ss} er seine Lehren angefochten sehen mu{\ss}te
(vgl.  Brief an die Kowalewska [sic] vom 24. M{\"a}rz 1885, mitgeteilt
von Mittag-Leffler, Acta Math. Bd. 39, S. 194ff).'' (Ibid.)}
\end{enumerate}
In this connection, Renate Tobies writes:
\begin{enumerate}\item[]
Klein and Lie discovered that they had many interests in common.  Like
Klein, Lie, too, shunned Weierstrass' lectures.%
\footnote{In the original: ``Klein und Lie entdeckten viele gemeinsame
  Interessen.  Wie Klein mied auch Lie die Vorlesungen von
  Weierstra\ss'' Tobies \cite[p.\;22]{To81}.}
\end{enumerate}
In 1919 Klein and Hilbert habilitated Emmy Noether, another central
figure in the development of modern mathematics, who strangely is only
briefly mentioned by Gray and not at all by HM.%
\footnote{In the foreword to her \emph{Habilitation} thesis that
contains the famous Noether's theorems, she writes that ``Klein's
second note and the present developments have been mutually influenced
by each other'' \cite{No71}.}

Ernst Zermelo, one of HM's \emph{moderns}, was habilitated by Hilbert
and Klein in 1899.%
\footnote{Ebbinghaus \cite[p.\;32]{Eb15}.}
It was at G\"ottingen, not at Berlin where he had studied, that
Zermelo first heard about set theory in a lecture given by Schoenflies
in 1898.  

Thus Klein appreciated the work of, and provided support for, numerous
modern mathematicians including Fraenkel, Lie, Noether, and Zermelo --
in some cases, against the opposition of the Berliners.

\subsection{Klein and Hilbert at G\"ottingen}
\label{s63}

Klein and Hilbert joined forces in turning G\"ottingen into a
broad-based, leading mathematical center offering a versatile
alternative to what was becoming a stifling influence of Berlin's
focus on arithmetic analysis.  Similarly, they jointly transformed the
\emph{Mathematische Annalen} into a leading research venue.  Contrary
to Quinn's claim (see Section~\ref{s4}), the emerging set-theoretic
foundations of the ``late 1800s'' found an avid supporter in the
person of Felix Klein:
\begin{enumerate}\item[]
One of the keys to th[e] success [of \emph{Mathematische Annalen}] was
Klein's knack for promoting the works of mathematicians who either
were estranged from or stood outside the mainstream influence of the
Berlin school: two prominent examples being Georg Cantor and David
Hilbert.%
\footnote{Rowe \cite[p.\;192]{Ro89}.  See further in \cite{Ro18}.}
\end{enumerate}
Klein published no fewer than 15 of Cantor's papers in
\emph{Mathematische Annalen}.%
\footnote{See Tobies \cite[p.\;38]{To90}.}
HM acknowledges that Klein also used Cantor as a referee for the
journal.%
\footnote{Mehrtens \cite[p.\;207]{Me90}.}
In fact, some of the mathematicians that would surely be listed as
\emph{modern} by HM, Gray, and Quinn turn out upon closer inspection
to be far more retrograde than Felix Klein on certain important
mathematical issues.  While the essence of mathematics may be its
freedom according to Cantor so that all mathematical theories are in
principle equally free, it turned out that some mathematical theories
were more equally free than others:
\begin{enumerate}\item[]
[Cantor] once called infinitesimals the ``cholera bacillus of
mathematics" and on another occasion referred to them as ``paper
magnitudes that have no other existence other than to be on the paper
of their discoverers and disciples."  Hilbert and Klein, on the other
hand, both recognized the validity of non-Archimedean systems in
geometry\ldots{}%
\footnote{Rowe \cite[pp.\;195--196]{Ro94}.}
\end{enumerate}
In this sense Felix Klein was arguably more modern than Georg Cantor.%
\footnote{See further in Kanovei et al.~\cite{18i}.}

If reading Quinn leaves one with the impression that Hilbert didn't
share Klein's sentiment concerning the importance of intuition, one
may be pleasantly surprised to read the following comment by Hilbert:
\begin{enumerate}\item[]
In mathematics, as in all scientific research, we encounter two
tendencies: the tendency toward abstraction - which seeks to extract
the logical elements from diverse material and bring this together
systematically - and the other tendency toward \emph{intuition}
[Anschaulichkeit], that begins instead with the lively comprehension
of objects and their substantial [inhaltliche] interrelations
{\ldots}\;The many-sidedness of geometry and its connections with the
most diverse branches of mathematics enable us in this way [namely,
  through the anschauliche approach] to gain an overview of
mathematics in its entirety and an impression of the abundance of its
problems and the rich thought they contain.%
\footnote{Hilbert as cited by Rowe in \cite[p.\;199]{Ro89}; emphasis
  added.}
\end{enumerate}
Hilbert valued intuition and its role in geometry, and wanted to
maintain that tradition in the 1920s when he taught his course on
\emph{Anschauliche Geometrie} four times before it was published with
Cohn-Vossen in 1932 as \cite{HC}.  Writes Rowe:
\begin{enumerate}\item[]
[Klein and Hilbert] were both elitists who were not afraid to talk to
the educated layman in simple terms {\ldots}\;It just was not
fashionable for German scholars to ``stoop to the masses''--indeed it
was regarded as a betrayal of the ideals of \emph{Wissenschaft}.
Klein paid a heavy price for doing so, having to suffer the scourn
[sic] of the famous mathematicians of the Berlin school--Weierstrass,
Kronecker, and Frobenius--who regarded him as a \emph{charlatan}.%
\footnote{Rowe \cite[p.\;202]{Ro94}; emphasis added.  Helmholtz,
citing Kronecker, explicitly referred to Klein as a ``charlatan'' and
none of the luminaries present, including Weierstrass, seem to have
objected (see Section~\ref{s44}).}
\end{enumerate}
The unjustified scorn for Klein felt by Weierstrass, Kronecker, and
Frobenius may have ultimately been an influence on the attitudes
displayed by a later generation of commentators such as HM, Gray, and
Quinn.

\section{Wrong address for countermodernism}

HM thought that he found a countermodern at G\"ottingen, in the person
of Felix Klein.  We argue that HM would have been more successful in
his search for countermoderns at another address.

\subsection{A meeting in Weierstrass' apartment in Berlin}
\label{s44}

On 14 january 1892, Adolf Hurwitz addressed a letter to Klein.
Kronecker having died in 1891, a position became available at Berlin.
Hurwitz wrote:
\begin{enumerate}\item[]
you can imagine how deep an impression the passing of Kronecker has
made on us here [in K{\"o}nigsberg].
%
%
Just about anyone in the \emph{unprejudiced} mathematical world agrees
that only you can fill the gap that has arisen in Berlin.  The only
uncertain [step in the process of appointing Klein to Berlin] is
whether [the decision-makers] inside the ministy are able to withstand
the countercurrents.%
\footnote{Hurwitz to Klein, cited in Tobies \cite{To99}.}
\end{enumerate}
Hurwitz appears to have expected the ministry to support Klein, and
anticipated resistance from ``countercurrents'' among mathematicians
at Berlin.

Eight days later, Berlin attitudes toward Klein were clearly on
display during a discussion that took place among leading Berlin
faculty members (the meeting took place in Weierstrass' apartment due
to his poor health).  Writes Rowe:%
\footnote{Rowe \cite[p.\;433]{Ro86}.  Another translation appears in
  Grattan-Guinness \cite[p.\;9, note\;22]{Gr09}.}
\begin{enumerate}\item[]
On 22 January [1892] a committee of Berlin faculty members met to
propose candidates to fill the vacancies, and on one point there was
unanimity of opinion.%
\footnote{Rowe's comment concerning alleged unanimity is not entirely
accurate as August Kundt, tenured physicist at the time of the
meeting, had spoken favorably of Klein during the meeting.}
As excerpts from the committee's protocol indicate, under no
circumstances would any of them countenance the candidacy of Felix
Klein:%
\footnote{Rowe's summary that follows is a shortened version of a much
longer transcript of the exchange that took place in 1892; see
\cite{Ts} and \cite[p.\;305--307]{Bi88}.  The order of the comments by
Helmholz and Weierstrass is reversed by Rowe.}
\vskip4pt
\begin{quote}
\textbf{Helmholtz}: Kronecker spoke very disparagingly of Klein.  He
regarded him as a charlatan [\emph{faiseur}].%
\footnote{The term \emph{faiseur} rendered by Rowe as ``charlatan'' is
a French loan word referring to someone involved in wheeling and
dealing, close in meaning to \emph{forger}, \emph{gambler}, or
\emph{cheat}.}
\vskip4pt

\textbf{Weierstrass}: Klein dabbles more. A bluffer [\emph{Blender}].%
\footnote{\label{f32}A fuller citation of the original is reproduced
in \cite[p.\;4]{Ts}: ``Schwarz bleibt bei der Stange.  Guter Vortrag.
Klein nascht mehr.  Blender.''  The word \emph{nascht} rendered by
Rowe as ``dabbles'' more literally means ``nibbles.''  The word
\emph{mehr} rendered by Rowe as ``more'' actually means ``rather'' in
context, as in ``Klein dabbles, rather.''  The word \emph{Blender}
rendered by Rowe as ``bluffer'' more literally means ``one who
dazzles'' with a pejorative implication of having little substance;
see Bair et al.~\cite{18b} for a similar use of the term
\emph{Schwindel} by Weierstrass in reference to what he described as
``the acolyte society'' of Klein, Lie, and Mayer.}
\vskip4pt

\textbf{Fuchs}: I have nothing against his person, only his pernicious
manner when it comes to scientific questions.
\end{quote}
\end{enumerate}
Whether or not Klein was a \emph{faiseur} and/or \emph{Blender}, the
tone of the 1892 discussion at Berlin is refreshingly frank, and
speaks volumes as to who exactly displayed \emph{countermodern}
tendencies at the time.

\subsection{Apartment meeting following Biermann}
\label{s62}

The meeting already discussed in Section~\ref{s44} provides evidence
of countermodern tendencies at Berlin.  The discussion of the
candidacy of Klein mentions a book of his (\emph{Lectures on the
Icosahedron}, 1884) but fails to mention his even earlier
\emph{Erlangen program} dating from 1872.  As noted in \cite{18b}, the
modernist nature of the \emph{Erlangen program} is recognized by most
commentators -- with the exception of HM.

The comments we cited in Section~\ref{s44} (following Rowe) suggest
that the participants in the apartment meeting were uncomfortable with
Klein's approach to mathematics and perhaps specifically with the
modernity of the \emph{Erlangen program}, which was not the
traditional kind of mathematics they were accustomed to.  Perusing the
discussion, it is hard to believe that the candidate Felix Klein
mentioned by Weierstrass at the outset:
\begin{enumerate}\item[]
at most Klein and Schwarz are eligible.  And [we] should only refrain
from trying [to get them to Berlin] if they are absolutely
unavailable.  Other personalities are not to be considered%
\footnote{Biermann \cite[p.\;305]{Bi88}; translation ours.}
\end{enumerate}
is the same candidate Felix Klein about whom the ensuing report to the
Prussian ministry said:
\begin{enumerate}\item[]
[F]irst and foremost it had to be considered [with foresight] that the
candidate shall be suited to continue what has been practiced at our
university for [several] generations: to lead the students to
\emph{serious, self-effacing, and in-depth} work on mathematical
problems. {\ldots}%
\footnote{Biermann \cite[pp.\;307--308]{Bi88}; translation ours;
  emphasis added.}
\end{enumerate}
The principles announced thus far certainly sound reasonable.  Their
application to the evaluation of Klein's candidacy may surprise the
modern reader:
\begin{enumerate}\item[]
{\ldots} For this reason [we] had to discard such personalities as
Professor Felix Klein in G\"ottingen (born 1849), on whose scientific
achievements the opinion of other scholars is very divided, but whose
actions as a whole, both through his writings and his teaching,
\emph{contradict the tradition of our university} just described.%
\footnote{Ibid.; emphasis added.}
\end{enumerate}
At the outset, Weierstrass describes Felix Klein, along with Hermann
Schwarz, as the two leaders of German mathematics who are the only
ones in the running for Kronecker's prestigious chair.  By contrast,
the final report presents Klein as a character of questionable record
both in research and in teaching who cannot even be relied upon to
supervise students in ``the tradition of the university'', being
apparently \emph{un}-serious, \emph{not} self-effacing
(self-centered?), and \emph{not} in-depth (shallow?).

What emerged at the meeting is that the Berlin traditionalists did not
trust Klein with the direction his mathematics (and possibly the
\emph{Erlangen Program} specifically) was taking, thereby
``contradicting the tradition'' of their university.  Signing that
report, the traditionalists in effect signed away the leadership of
the Berlin school in mathematics for decades to come.

\section{Geometry, intuition, logic, physics}
\label{s33b}

Klein and Hilbert shared a belief in a decisive role of intuition in
mathematics.  Thus, Hilbert believed his \emph{Grundlagen der
  Geometrie}~\cite{Hi99} to be a study of our geometric intuition.

\subsection{Hilbert's \emph{Grundlagen}}

In the Preface to his \emph{Grundlagen}, David Hilbert writes:
\begin{enumerate}\item[]
The choice of the axioms [of geometry] and the investigations of their
relations to one another is a problem which, since the time of Euclid,
has been discussed in numerous excellent memoirs to be found in the
mathematical literature.  \emph{This problem is tantamount to the
  logical analysis of our intuition of space}.  (emphasis added)
\end{enumerate}
In Klein's remarks on Hilbert's \emph{Grundlagen}, he stresses the
importance of logic in geometrical investigations:
\begin{enumerate}\item[]
I have already emphasised one thing about which most people today are
in reasonable agreement.  That is that we are concerned here with the
leading concepts and statements, which one must of necessity put into
the front rank of geometry in order to be able to realise mathematical
proofs derived from them by pure logic.  This statement does not
answer the question as to the real source of these leading concepts
and theorems.  There is the old point of view that they are the
intuitive possession of every person, and that they are of such
obvious simplicity that no one could question them.  This view,
however, was shaken, in large measure, by the discovery of
non-Euclidean geometry; for here it is clearly shown that space
intuition and logic by no means lead compellingly to the Euclidean
parallel axiom.%
\footnote{Klein \cite[p.\;213]{Kl2016}.}
\end{enumerate}
Having thus analyzed the relation of logic and intuition, Klein
continues:
\begin{enumerate}\item[]
To the contrary, we saw that, with an assumption, which contradicts
the parallel axiom, we come to a logically closed geometric system,
which represents actual perceptual relations with any desired degree
of approximation.  However, it may well be claimed that this parallel
axiom is the assumption, which permits the simplest representation of
space relations.  Thus it is true in general that fundamental concepts
and axioms are not immediately facts of intuition, but are
appropriately selected idealisations of these facts.  The precise
notion of a point, for example, does not exist in our immediate
sensory intuition, but is only a fictitious limit, which, with our
mental pictures of a small bit of shrinking space, we can approach
without ever reaching.%
\footnote{Ibid.}
\end{enumerate}
Thus Klein, like Hilbert, stressed the role of logic in the study of
mathematical intuition.

\subsection{Physics as seen by Hilbert and Klein}

Mehrtens, Gray, and Quinn sought to portray Klein as looking backward
to the 19th century.  Their critique imputes a further allegedly
\emph{countermodern} shortcoming to Klein, namely an interest in
applications to physics.  This debate (concerning a proper relation of
mathematics to physics) we suspect will never end, and it was the real
issue at the root of the debate sparked by the text by Jaffe and Quinn
\cite{JQ}, leading some physicists to coin the term ``rigorous
mathematics in the sense of Jaffe and Quinn.''%
\footnote{See e.g., Eager and Hahner \cite{Ea23}.}

Meanwhile, anyone who seeks to recruit Hilbert to a hypothetical
\emph{modernist} cause will have to come to grips with Hilbert's
profound interest in physics in general and Einstein's relativity
theory, in particular.%
\footnote{See \cite{18b}.}
As acknowledged by Gray, 
\begin{enumerate}\item[]
[Einstein] found he got a better hearing from Hilbert and Klein in
G\"ottingen than he did from his colleagues in Berlin.%
\footnote{Gray \cite[p.\;326]{Gr08}.}
\end{enumerate}
Concerning the issue of the proper relation between mathematics and
physics, one may well ponder the following story told by one of the
world's leading mathematicians in 2014:
\begin{enumerate}\item[]
In a mythical country, there was a very skilled blacksmith who could
make extraordinary swords, and extraordinary axes. (I visited the
weapon museum in Nagoya, and I can tell you that the Japanese
blacksmiths were very skilled indeed.)  One day, the blacksmith's son
used an axe from his father's collection to cut down a dangerous tree.
The father was furious, and told him that this axe was not to be used,
for it was a museum piece.  The tale is easy to understand.%
\footnote{Quoted in Paycha \cite{Pa17}.}
\end{enumerate}
The speaker is none other than Pierre Cartier.

\subsection{Differences according to Max Born}
\label{s57}

Physicist Max Born, one of the fathers of quantum mechanics, was
trained at G\"ottingen.  Born commented as follows on the differences
in style between Klein and Hilbert:

\begin{enumerate}\item[]
[Klein] was the most dazzling lecturer, but he was too brilliant for
me.  One would have to have been a really pure mathematician to enjoy
his lectures, even when he treated physical or technological
applications, which he liked to do. \ldots{} At this time of my life, I
was more inclined to pure mathematical rigor - inclined to
epsilontics, as we called it in our own jargon. \ldots{} $\varepsilon$
almost never appeared in Klein's lectures.

If Klein arrived at such a critical point, he used to say: 'I will
sketch for you the idea of the matter; the rigorous proof can be
looked up in one of the following publications', and he then gave a
literature reference.  He even added a disdainful remark about those
people who need a proof of such an obvious fact.  


Now, the proofs that Hilbert gave in his lectures were rather sketchy,
too, but I didn't ever have any difficulties to reconstruct them
according to his clear indications.

In the case of Klein, however, I often got hopelessly lost 
%
%
and was force to read a large number of boring treatises. 
%
%
In their general structure, Hilbert's and Klein's lectures were very
different.

Hilbert was a mountain guide, who led one on the shortest and most
secure route to the summit; Klein, on the contrary, was rather like a
prince, who wanted to show his admirers the size of his kingdom, by
leading them on endlessly winding paths, and who stopped on every
small hill to give an overview of the ground already covered.  I got
impatient with this method [of lecturing] and wanted to reach the
summit quickly.  Therefore I preferred Hilbert.%
\footnote{Born \cite[p.\;131]{Bo75}.}
\end{enumerate}
Born then reveals his punchline:
\begin{enumerate}\item[]
Nowadays, my taste has changed.  Whenever I have the time to read a
mathematical work purely for my own pleasure, I reach for one of
Klein's books \ldots{} and I take delight in them like [one takes
delight] in a work of art, a great novel or a biography. Indeed, some
of these books of Klein's are an integral part of the culture in which
I grew up, and which is now going to pieces.%
\footnote{Ibid.}
\end{enumerate}
Felix Klein's style as described by Born is similar to that of Mikhael
Gromov, one of the greatest and most influential geometers living.%
\footnote{For example, Gromov's book \cite{Gr} has been cited by over
2800 articles; see
\url{https://scholar.google.com/scholar?cites=2060554209746836258}}
%
%
%
Clearly, the differences in style between Klein and Hilbert persist in
modern mathematics.  The attempt by Mehrtens, Gray, and Quinn to
characterize one of them as countermodern is therefore unconvincing.

\subsection{Hilbert's formalism}
\label{s6}

Corry warns against simplistic interpretations of Hilbert's views on
formalism:
\begin{enumerate}\item[]
David Hilbert is widely acknowledged as the father of the modern
axiomatic approach in mathematics.  The methodology and point of view
put forward in his epoch-making \emph{Foundations of Geometry} (1899)
had lasting influences on research and education throughout the
twentieth century.  Nevertheless, his own conception of the role of
axiomatic thinking in mathematics and in science in general was
significantly different from the way in which it came to be understood
and practiced by mathematicians of the following generations,
including some who believed they were developing Hilbert's original
line of thought.%
\footnote{Corry \cite[Abstract]{Co07}.}
\end{enumerate}
Furthermore,
\begin{enumerate}\item[]
  Among the important sources of ideas that inspired the New Math%
\footnote{For Freudenthal's critique of New Math, Piaget, and Bourbaki
see Freudenthal \cite{Fr91} and Katz \cite{20e}.}
was a certain perception of the Moore Method and the attempt to apply
to school mathematics what this method had considered to be of high
value in the training of research mathematicians.%
\footnote{Ibid., Section~1.}
\end{enumerate}
The difference between Hilbert's conception and that of some
mathematicians of the following generations who ``believed they were
developing Hilbert's original line of thought'' is well illustrated by
the case of Frank Quinn.  While championing what he sees as Hilbert's
formalist approach, Quinn reproaches Hilbert for having allegedly
\begin{enumerate}\item[]
accepted the slanders by saying `mathematics is a game played
according to certain rules with meaningless marks on paper.'%
\footnote{Quinn \cite[p.\;35]{Qu12}.}
\end{enumerate}
Quinn goes on to venture an alternative line of argument to Hilbert's
that he feels would have been historically more effective.%
\footnote{The text of such an alternative runs as follows: ``Axiomatic
  definitions can be artificial and useless, but they can also
  encapsulate years, if not centuries, of difficult experience, and
  newcomers can extract reliable and effective intuitions from
  them. Similarly, fully detailed arguments can be formal and
  content-free, but fully confronting all details usually deepens
  understanding and often leads to new ideas. Fully detailed arguments
  also give fully reliable conclusions, and full reliability is
  essential for successful use of the powerful but fragile
  excluded-middle method'' \cite[p.\;35]{Qu12}.  We will withhold
  judgment as to the effectiveness of Quinn's proposal.}
However, Quinn's claim attaches exaggerated importance to Hilbert's
comment.  Hilbert may have gotten carried away in the heat of the
argument against Brouwer to make the comment on \emph{meaningless
  marks}.  Hilbert's comment does not represent the thrust of
Hilbert's position, as argued by Avigad and Reck~\cite{AR}.  Hilbert's
comment should be understood in the context of his
\emph{metamathematical} program rather than as a statement of a
\emph{philosophical} position.  The alleged philosophical significance
of the comment has been exaggerated out of proportion subsequently, by
many authors besides HM and Quinn.  In defense of Hilbert,
Novikov~\cite{No2} argues that Hilbert was as dedicated to
\emph{meaning} as any great mathematician.%
\footnote{See \cite{11a} for a related discussion.}
Meanwhile, Gray declares:
\begin{enumerate}\item[]
By holding onto meaning, [Kronecker and Enriques] belong to the group
of mathematicians Mehrtens called the \emph{Gegenmoderne}.%
\footnote{Gray \cite[p.\;274]{Gr08}.}
\end{enumerate}
Accordingly, ``holding on to meaning'' constitutes a shortcoming and a
litmus test of being a \emph{Gegenmoderne}.  But by Gray's criterion,
Hilbert would also be among the \emph{Gegenmoderne}, as he arguably
adhered to meaning as much as Kronecker, Enriques, and Klein.  Pinning
modernism to an absence of, or indifference to, \emph{meaning} is a
position that is too simplistic to be useful as a historiographic
tool.

\section{Conclusion}
\label{s8}

We have argued that Klein and Hilbert were \emph{modernist} allies in
a common cause challenging what in retrospect appear to be artificial
limitations on mathematical research imposed by some of the
mathematical trend-setters at Berlin, which would have been a more
appropriate address than G\"ottingen for a search for countermoderns.
Klein's and Hilbert's challenge found a clear expression in Hilbert's
famous list of 23 problems, only a minority of which were concerned
with arithmetized analysis as emphasized by the Berlin school.

Mehrtens' modern/countermodern dichotomy is artificial and appears to
be tailor-made to cast mathematicians like Klein in unfavorable light.
Mehrtens' distaste for Klein stems from Mehrtens' political
affiliations (see Section~\ref{s21} for Rowe's comments on Mehrtens'
dialectics) rather than his historical scholarship.  The cause of
Mehrtens' resentment of Klein appears to be Klein's role in the
academic-industrial complex of Wilhelmine Germany of the pre-world war
1 period.

Mehrtens' decision to incorporate the interest in applied mathematics
as a characteristic of the so-called \emph{countermoderns} is a thinly
veiled attempt to design a theoretical framework that would make Klein
look retrograde.  Mehrtens' dialectical methodology influenced authors
ranging from Gray and Quinn to, more recently, Siegmund-Schultze and
Mazzotti.  Not all of these authors share Mehrtens' original political
motivation.  Thus, Siegmund-Schultze reports ``shock'' as his first
reaction upon reading an early draft of Mehrtens' book (see
Section~\ref{s21b}).  With Rowe and Tobies, we reject Mehrtens'
provocative linkage of Klein's ethnographic ideas to Vahlen-style
antisemitism, as well as Gray's and Mazzotti's toleration of such
linkage.

We have examined the attempts by Mehrtens, Gray, Quinn, and to a
certain extent Siegmund-Schultze to question Klein's modernity.  We
concur with Grattan-Guinness, Rowe, and Tobies that the portrayal of
Klein and Hilbert as opponents in the matter of modernity is contrary
to much historical evidence.

\section{Acknowledgments}

We are grateful to Hollis Williams for helpful criticisms of an
earlier draft of the article.

\end{document}